\documentclass[article,12pt]{elsarticle}
\usepackage{amssymb}
\usepackage{mathrsfs}
\usepackage{exscale}
\usepackage{makeidx}
\usepackage{graphicx}
\usepackage[T1]{fontenc}
\usepackage{amsmath,amssymb}
\usepackage[all]{xy}

\newtheorem{thm}{Theorem}
\newtheorem{rem}{Remark}
\newtheorem{df}{Definition}

\newtheorem{lem}{Lemma}
\newtheorem{exemple}{Example}

\newcommand{\preuve}{\indent {\it Proof:}\hspace{4mm}}

\journal{Annals of pure and applied logic}

\begin{document}

\begin{frontmatter}



\title{Algebraically 
closed structures in Positive Logic.}


\author{Mohammed Belkasmi}
\ead{belkasmi@math.univ-lyon1.fr}
\address{Maths Department\\
 Qassim University\\
 PO Box 6644 Buraidah 51452\\
Saudi Arabia }

\begin{abstract}

In this paper we  extend  
of the notion of  algebraically 
closed given in the case of groups and skew fields to 
an arbitrary h-inductive theory. The main subject of this paper  is the study 
of the notion of positive algebraic
closedness  and its relationship
with the notion of positive closedness 
and the amalgamation property.
\end{abstract}

\begin{keyword}
Positive logic, h-inductive theory, h-universal theory
positively existentially closed, positively  
closed , positively algebraically 
closed, e-elementary extension. 

2010 Mathematics Subject Classification: Primary 03C95;
Secondary 03C48, 03C10

\end{keyword}

\end{frontmatter}




\section{Positive logic}
\subsection{Introduction}
The main technical tools developed by Abraham Robinson in
 the study of the notions 
of model-complete theories, existentially closed models  and  
inductive theories in the framework of first 
order logic with negation,  are the notions of 
inductive
sentence, embedding,  existential formula and universal sentences
which are precisely negations of existential sentences. 

The positive logic provides a simpler framework than that of  
the logic with negation, in the sense that only positive existential
formulas and  h-inductive sentences are considered.  
The inductive theories in Robinson's sense 
 become 
a  special positive case,    by a partial
morleysation adding a new relation symbol for the negation of each
atomic formula of the language.

Positive logic offers really new situations; it permits
the
extension to homomorphisms the  model-theoric notions which are classically
associated to embeddings, in particular the notion of existentially closed models.
It permits also the restriction of the set of formulas to the 
positive ones.

The present paper provides an  exposition of the positivisation
of the notions of existentially and algebraically closed structures.
This approach turns out to be compatible with some algebraic notions, and create new opportunities.  
 
The plan of this  paper is as follows. In the section 1, 
after revising 
the foundation of the positive logic, we give a generalized form 
of amalgamation called 
the amalgamation embedding-immersion, which turns out to be very useful in the  study the properties of the notions introduced in the section 2.\\
 The section 2 is divided into two parts. In the first part, 
 we introduce the notion of positively  algebraically closed structure, we 
 analyse the preservation of this notion by various type of extensions, and by passing to substructures.\\
Afterwards we give a 
  syntactic characterisations of 
 positively algebraically closed structures
  and the theories with algebraic model-companion.
  we finish this part by analysing  
  the relation between the classes of positively  algebraically closed 
  and existentially closed structures.
  In the second part, we define the 
 notion of e-elementary algebraic extension and prove some 
 fundamental properties of this notion. 
 
 \subsection{Pc models}
Our terminology and notations will be consistent
with \cite{ana} and \cite{begnacpoizat}.\\
Let $L$ be a first order language.
The positive formulas are formulas which are obtained from the
atomic formulas by the use of\  $\wedge, \vee $ and $\exists$.
They are of the form 
 $\exists\, \bar x\,\psi(\bar x,\bar y)$, where 
$\psi$ is  positive quantifier-free; the variables  
$\bar y$ are said to be free.

 A sentence is said to be h-inductive 
 if it is a finite conjunction of sentences 
of the form: $$\forall\bar x\ \ \exists
\bar y\psi(\bar x, \bar y)\rightarrow\ \ 
\exists\bar z\varphi(\bar x,\bar z)$$

where $\psi, \varphi$  are positive
 quantifier-free formulas.

The h-universal sentences represent a  special case of
h-inductive sentences; they are the sentences that can be written
as  negation of a positive sentence.
 
Let $A$ and $B$ be two $L$-structures. By an homomorphism $h$
from $A$ to $B$ we mean a mapping $h$ from $A$ to $B$  such that,
for every  tuple $\bar a$ from  $A$ ($\bar a\in A$ by abuse of notation)
and for every atomic formula $\phi$. $A\models\phi(\bar a)$
implies $B\models\phi(h(\bar a))$.
In such a case, $B$ is said to be a continuation of\  $A$.\\
By an embedding of\  $A$ into $B$ we mean an homomorphism 
$h:\ A\longrightarrow B$ such that  for every
atomic formula $\phi$, $A\models\phi(\bar a)$
if and only if $B\models\phi(h(\bar a))$. The homomorphism 
$h:\ A\longrightarrow B$
 is an immersion whenever for every  $\bar a\in A$
and $h(\bar a)$ satisfy the same $L$-positive formulas.

A class of $L$-structures is said to be $h$-inductive 
(resp. $e$-inductive, $i$-inductive) if it is closed
under inductive limits of homomorphisms
(resp. embeddings, immersions). It is
easy to verify that the class of models
of an $h$-inductive theory is $h$-inductive.
The theorem 23 of \cite{begnacpoizat}
shown that this is indeed a characterization of 
the $h$-inductive classes of models of a first order theory.

The positively closed  structures  are 
the central objects of study in positive model theory. 
The corresponding notion in the model theory with negation 
are the existentially complete structures.
Broadly speaking, the positively closed structures  are the structures for which the h-universal sentences
 are persistent under homomorphisms in the sense that, 
if a positive formula is satisfied in any continuation of 
the structure, then it is satisfied in the structure itself.
\begin{df}
Let $L$ be a first-order language and $\Gamma$ a class  of\  $L$-structures.
A member $M$ of\  $\Gamma$
is said to be positively closed  (pc from
now on),
if every homomorphism from $M$ into an element
of\  $\Gamma$ is an immersion.
\end{df}

The following facts  give  fundamental properties of 
pc models. They will
 be used  without mention in the rest of this paper.
\begin{lem}( Th\' eor\` ereme 2, {\cite{begnacpoizat}})\label{pecconti}
Every member of an $h$-inductive class 
has a pc continuation in the same class.
\end{lem}
\begin{lem}(Lemme 14, \cite{begnacpoizat})
A model $B$ of an h-inductive $L$-theory $T$ is a pc model of 
$T$ if and only if
for every $\bar a\in A$ and $\varphi$ a positive $L$-formula, if 
$A\nvDash\varphi(\bar a)$ then there is $\psi$ a positive 
$L$-formula such that,
$$
 \left\{
    \begin{array}{ll}
         A\models\psi(\bar a)\\
       T\vdash \neg\exists\bar x\, (\varphi(\bar x)\wedge\psi(\bar x)).
    \end{array}
\right.
$$
\end{lem}

For every  formula $\phi$, we denote by   
$Ctr_T(\phi)$ the  set of  positive formulas $\psi$ such that 
$T\vdash\neg\exists \bar x(\phi(\bar x)\wedge\psi(\bar x))$.\\
Given $A$ be a $L$-structure, we let 
 $L(A)$ be the language obtained from $L$ by adjoining 
the element of\  $A$ as constants. 
We denote by:
\begin{itemize}
\item $Diag(A)$ the 
set of   atomic and negated atomic $L(A)$-sentences  satisfied by $A$.
\item $Diag^+(A)$ the set of positive quantifier-free $L(A)$-sentences 
satisfied by $A$.
\end{itemize} 
\begin{df}
\begin{itemize}
\item Two h-inductive theories are said to be companion 
 if they 
have the same pc models.
\item An h-inductive theory $T$ is said to be positively complete
(or it has 
the joint continuation property) if  
 any two models of T have a common continuation.
\end{itemize}
\end{df}
Every $h$-inductive  
theory $T$ has a maximal companion denoted $T_k(T)$,  called
the Kaiser's hull of\  $T$.  $T_k(T)$
is the set of  $h$-inductive consequences of the pc models of\  $T$.
 Likewise,
$T$ has a minimal companion denoted $T_u(T)$, formed by its
$h$-universal consequences. 
\begin{rem}
 Let $T_1$ and $T_2$  be two h-inductive theories. 
The following  propositions are equivalent:
\begin{itemize}
\item $T_1$ and $T_2$ are companion.
\item $T_k(T_1)\equiv T_k(T_2)$.
\item $T_u(T_1)\equiv T_u(T_2)$.
\end{itemize}
\end{rem}

\begin{df}
Let $T$ be an h-inductive theory.
\begin{itemize}
 \item $T$ is said to be model-complete 
if  every model of\  $T$ is a pc model of\  $T$.
\item We say that $T$ has a model-companion whenever $T_k(T)$ is 
model-complete.
\item An $n$-type is a maximal set of positive formulas
in $n$ variables that is consistent with $T$. We denote by $S_n(T)$  the space
of\  $n$-types of the theory $T$. 
\end{itemize}
 \end{df}

\begin{exemple}\label{exemplepecs}
Let $L$ be the language consisting of one unary function $f$.
\begin{enumerate}

\item Let $T$ be the  h-inductive $L$-theory 
$$\{\forall x, y\, (f(x)=f(y)\longrightarrow\, x=y)\}.$$ 
Every model of  $T$ 
can be viewed as directed graph formed by   union of cycles,
  two-sided and one sided chains. The 
 elements of the model  are  the vertices of the  graph, and 
 two elements $a$ and $b$ of the model are
connected by a directed edge pointing from $a$ to $b$, if the model 
satisfies the formula $f(a) = b$.\\
$T$ has a model-companion, and the class of pc models of\  $T$
 is  reduced to  the structure $A_e=\{x\}$. 
\item Let 
$T'$ be the h-universal $L$-theory  
$$\{\neg\exists x\, f\circ f(x)=x\}.$$
Every  model  $(A, f)$ of\  $T'$ can be represented by a  
directed graph as in the previous example.
 $T'$ has a unique 
 pc model formed by one  cycle of
order 4 and one cycle of order p for every prime $p > 2$. Thereby the class of pc models is not elementary.
\item This example is given in  \cite{Almaz}. Let $T''$ be the 
h-universal theory $\{\neg\exists x\, f(x)=x\}$. 
$T''$ has a unique 
 pc model formed by a cycle of order p for every prime $p$.
 \item Let $T_G$ be the theory of groups in the usual
 language of groups. The class of  pc
model of\  $T_G$ is  reduced to the trivial group.
\end{enumerate} 
\end{exemple}
Let $A$ be a L-structure. We shall use the following notations:
\begin{itemize}
\item $T_i(A)$ (resp. $T_u(A)$) denote the set of h-inductive
 (resp. h-universal ) 
 $L(A)$-sentences satisfied by $A$.
\item $T_i^\star(A)$ (resp. $T_u^\star(A)$) 
denote the set of h-inductive (resp. h-universal) $L$-sentences 
 satisfied by $A$.
 \item $T_k(A)$ (resp. $T_k^\star(A)$) denote the Kaiser's
hull of\  $T_i(A)$
(resp. of  $T_i^\star(A)$).
\end{itemize}
 \begin{rem}
 Let $A$ and $B$ be  two $L$-structures such that $A$ is 
 immersed in $B$, then $T_i^\star(B)\subseteq T_i^\star(A)$.
 \end{rem}
 \begin{lem}
 Let $B$ a model an h-inductive.
  If $B$  is a continuation of a pc model $A$ of\  $T$, 
 then $A$ is a pc model of 
 $T_i^\star(B)$.
 \end{lem} 
 \preuve 
 Since $A$ is immersed in $B$ then $A$ is a model of 
  $T_i^\star(B)$. Given that $T\subseteq T_i^\star(B)$
  it follows that $A$ is a pc model of $T_i^\star(B)$.\qed
  
The following lemma provides a  generalized form of 
asymmetric amalgamations called embedding-immersion amalgamation, 
inspired by other forms of  asymmetric amalgamations
given in \citep{begnacpoizat} and \citep{ana}.
 This form of amalgamation turns out 
to be a very useful technical tool in the next section. 
\begin{lem}(embedding-immersion amalgamation)\label{pamalgam}\\
Let $A, B$ be two $L$-structures and 
$C$ a model of\  $T_i(A)$. Let
$e$ an embedding from $A$ into $B$ and  $i$  the natural immersion
from $A$ into $C$.
  Then there exists $D$ a model of\  $T_i(B)$
 such that the following
diagram commutes:
\[
\xymatrix{
    A \ar[r]^{e} \ar[d]_{i} & {B} \ar[d]^{i'} \\
    C \ar[r]_{e'} & {D}
  }
\]  
where $e'$ is an embedding and $i'$ is an immersion.
\end{lem}
\preuve
Consider the language $L^+$ formed by the language $L$ and 
the elements of\  $A, B$
and $C$, so that the elements of\  $A$ are 
interpreted by  the same 
symbol of constants in $B$ and $C$. The
proof consists to show that the  set of h-inductive $L^+$-sentences
   $T_i(B)\cup Diag(C)$ is consistent.
   
 Let $\Delta$ be a finite fragment of\  $Diag(C)$. 
 Without loss of generality, we may assume  that 
$\Delta=\{\neg\psi(\bar a, \bar c), \varphi(\bar a, \bar c)\}$
where $\psi$ and $\phi$ are positive quantifier-free formulas,
$\bar a\in A$ and  $\bar c\in C$. To show 
the consistency of  $T_i(B)\cup Diag(C)$, we will give an 
interpretation of\  $\Delta$ in $B$.\\
 We claim that 
  $A\nvdash\forall\bar y\, (\varphi(\bar a,\bar y)\rightarrow
\psi(\bar a, \bar y))$.  Indeed, if $A\vdash\forall\bar y\, (\varphi(\bar a,\bar y)\rightarrow
\psi(\bar a, \bar y))$, 
given that  $C\vdash T_i(A)$
 then  $C$  
satisfies $\forall\bar x\bar y\, (\varphi(\bar x,\bar y)\rightarrow
\psi(\bar x, \bar y))$, contradiction with 
$\{\neg\psi(\bar a, \bar b), \varphi(\bar a, \bar b)\}
\subset Diag(C)$.
Thereby there is $\bar a'\in A$ such that   
$A\nvDash\psi(\bar a, \bar a')$ 
and  $A\models\varphi(\bar a, \bar a')$.
Now  since  $A$ is embedded in $B$,  we can interpret $\varphi(\bar a, \bar a')$
and $\neg\psi(\bar a, \bar a')$ in $B$. \qed
\section{Positively Algebraically closed structures} 
Algebraically 
closed structures plays an important part in the study of
 algebraic theories in the framework of model theory with 
negation. The notion of positively  algebraically closed structures presented 
here is the outcome of the extension to an arbitrary h-inductive 
theory of the notion of algebraically closed groups introduced 
by Scott \cite{scott}, and considered in \cite{rob1} by A. Robinson
in the case of groups and skew-fields.

In the first part of this section, after having presented the notion of 
positively algebraically closed structures in adequate setting, we 
establish the fundamental properties of the class of 
positively algebraically closed structures of an h-inductive theory.

\begin{df}
Let $T$ be an h-inductive theory, a model $A$ of\  $T$ is said to be
positively  algebraically closed ( pac from now on) if every embedding 
from $A$ into a model of\  $T$ is an immersion.
\end{df}
\begin{rem}
\begin{itemize}
\item Every pc model of\  $T$ is a pac model of $T$,  
 thereby every model of\  $T$ is continued 
in some pac model of\  $T$.
\item  Let $T$ be an h-inductive theory such that, 
the negation of every  positive quantifier-free formula 
is equivalent modulo $T$ to a positive formula. Then 
every pac model of\  $T$ is a pc model,  and 
an existentially closed model of\  $T$.
\item  Let $L$ be a language without relation symbols. Let 
$T$ be an h-inductive theory such that the inequality is 
positively defined. Then the negation of every 
positive quantifier-free formula is positively defined. Thereby 
every ac model of\  $T$ is a pc model of\  $T$ and 
an existentially closed model of\  $T$.
\end{itemize}
 \end{rem}
The content of the notions of  pc and pac models 
are  best seen by considering the following  examples.
\begin{exemple} \label{exemplwpac}
\begin{enumerate}
\item Let $T$ be the h-inductive theory given in the example 
\ref{exemplepecs}.    
A model $(A, f)$ of\  $T$ is  pac 
 if and only if $f$ is bijective and has a fixed point.
 Thereby  the class of  pac models of\  $T$ is elementary.
 \item Let $T'$ be the h-inductive theory given in the example 
\ref{exemplepecs}.  A model of\  $T'$ is  pac if and only if 
it contains a $4$-cycle and a $p$-cycle for each prime $p>2$. 
Thereby the  pac models of\  $T'$  form an 
elementary class.\\ 
However, in the framework of 
model theory with negation, a model of\  $T'$ is an existentially 
closed model if and only if contains the set of\  $n$-cycles  
for each integer $n>2$. 
\item Every $L$-structure $A$ is a pac model of the theory $T_u(A)$
 over the language 
$L(A)$. 
\item Let $T_{ag}$ the theory of abelian groups
 in the usual language of groups. Every divisible abelian group 
 is a pac model of $T_{ag}$. However, 
 the class of existentially closed models of $T_{ag}$ is the 
 class of all divisible abelian groups that contain for each 
 prime number $p$ an infinite number of elements of order $p$. 
 
 \item Given that the language of the theory of fields is without 
 relation symbols, and the negation of the formula 
 $x=y$ is defined by the positive formula 
 $\exists z ((x-y).z=1)$. Then, every pac field  
  is   pc  and an existentially closed field. 
 \end{enumerate} 
 \end{exemple}
\begin{lem}
Let $T$ be an h-inductive theory. The class of 
pac models of\  $T$ is e-inductive (inductive under  embeddings).
\end{lem}
\preuve Let $A$ be the inductive limit of the inductive sequence
$(A_i, p_{i,j})_{i,j<\omega}$, where 
$A_i$ are
 pac models of  $T$ and  $p_{i,j}$   embeddings.
 Consider 
$B$ a model of\  $T$ and $p$ an embedding from $A$ into $B$. In order 
to show that $A$ is a pac model of\  $T$, 
 we will  show that $p$ is an immersion.\\ 
Suppose that 
$B\models\exists \bar x\,\psi(p(\bar a), \bar x)$ where $\psi$ 
is a positive quantifier-free
 formula, and $\bar a\in A$. Let $i<\omega$ such that 
$\bar a\in A_i$,  let $p_i$  be the natural
 homomorphism from $A_i$ into $A$.  Given that  $A_i$ is a pac model 
and $p\circ p_i$ is embedding then  
$p\circ p_i$ is an immersion. Since   
$p(\bar a)=p\circ p_i(\bar a)$, it 
follows that $A_i\models \exists \bar x\,\psi(\bar a, \bar x)$. Thereby  
$A\models \exists\bar x\,\psi(\bar a, \bar x)$, so $p$ is an immersion.\qed
\begin{lem}\label{embdweak}
Every model of an h-inductive theory $T$ is embedded in a  
pac model of\  $T$.
\end{lem}
\preuve
The proof is similar to the classic one for the 
existence of existentially closed models,
 except that the existential formulas must
  be replaced by the positive formulas.
  
 Let $A_0$ be a model of $T$. Let $\Gamma_0=\{\varphi_i\mid i<\alpha_0\}$ 
be the set of positive 
$L(A_0)$-formulas  enumerated   by an ordinal 
$\alpha_0$.
first, we  construct an inductive sequence $(M_i)_{ i<\alpha_0}$ of 
 models of $T$ as follows;
if the first formula $\varphi_0\in \Gamma_0$ is satisfied in some  model $B$
 of $T$ in 
which $A$ is embedded, we take $M_1=B$, if not we take $M_1=M_0=A_0$.
We continue in this manner,
if the second formula of $\Gamma_0$ is satisfied in some model 
$C$ of $T$ in which $M_1$ is embedded, we take $M_2=C$, if not
we take  $M_2=M_1$.  
 If
 $\beta<\alpha_0$ is a limit ordinal,  one defines  $M_\beta$
 as the 
inductive limit of $(M_i)_{i<\beta}$.\\
 Let $A_1$ be the inductive limit  of 
$(M_i)_{i<\alpha_0}$. 
We repeat the previous construction for $A_1$ and $\Gamma_1$ 
 the set of positive formulas with parameters in $A_1$. 
 We obtain an inductive
sequence  (under embeddings)  $(A_i)_{i<\omega}$.
 The inductive limit $C$ of  $(A_i)_{i<\omega}$ is a 
 pac model  of $T$, and $A$ is embedded in $C$.\qed

We use $A_T$ and $E_T$ to denote respectively 
the classes of pac and pc models of an h-inductive theory $T$. We have 
$E_T\subseteq A_{T_k(T)}\cap A_{T_u(T)}$.

 \begin{df}
Two h-inductive theories are said to be e-companion 
if  every model of one of them can be 
embedded into a model of
the other.
\end{df}  

  \begin{lem}
  Two h-inductive theories are e-companions if and only if 
  they have the same class of  pac models.
  \end{lem}
  \preuve Let $T_1$ and $T_2$ be two h-inductive theories with 
  the same class of pac models. 
  By the lemma \ref{embdweak}, every model of one of them can be 
embedded into a model of
the other.
  
  For the other direction,  assume that $T_1$ and $T_2$ are 
  e-companions theories.
  Let $A$ be a pac model of\  $T_1$.
   By the lemma \ref{embdweak}, we obtain  the 
 following diagram: 
 $$A\xrightarrow{e_1}B\xrightarrow{e_2}C\xrightarrow{e_3}D$$ 
 where $e_1, e_2, e_3$ are embeddings, $B$ a model of\  $T_2$,
 $C$ a pac
  model of
 $T_2$ and
  $D$  a model of\  $T_1$.
 Given that  $A$ is a pac model of\  $T_1$, then $e_3\circ e_2\circ e_1$ is 
 an immersion, so $e_1$ is an immersion. Consequently 
 $A$ is a model of\  $T_2$ and  a pac model
 of\  $T_2$. Likewise every  pac model of\  $T_2$ is a pac model of\  $T_1$.\qed
 \begin{lem}\label{substrextension}
 If $A$ is a substructure of a pac model $B$ of\  $T$ and 
    $B$ a model of\  $T_i^\star(A)$, then $A$ is a pac model 
 of\  $T$.
 \end{lem}
 \preuve
 Suppose that $B\models T_i^\star(A)$, let $p$ be an embedding from 
 $A$ into  a model $C$ of\  $T$. By  the lemma \ref{pamalgam}, we 
 have the following diagram: 
 \[
\xymatrix{
    A \ar[r]^{i} \ar[d]_{p} & {B} \ar[d]^{p'} \\
    C \ar[r]_{i'} & {D}
  }
\]
where $D$ is a model of\  $T$, $p'$ an embedding and $i'$ an immersion.
Given that $B$ is a pac model of\  $T$ then $p'$ is  an immersion, and 
so does $p$. Thereby $A$ is a pac model of\  $T$.\qed
 
The following theorem gives a syntactic characterisation of 
pac model,  based on a formal description of the positivation 
of negative formulas in a pac model.
 \begin{thm}\label{weakchara}
Let $A$ be a model of an h-inductive theory $T$. 
 $A$  is a  pac  model of\  $T$ if and only if 
for every $\bar a\in A$, and $\psi$  a positive formula such that 
$A\nvDash \psi(\bar a)$, there exist $\bar b\in A$ and 
two positive quantifier-free formulas 
$\theta_1, \theta_2$ such that:
\begin{enumerate}
\item $A\models 
\theta_1(\bar a, \bar b)$, and $A\nvDash\theta_2(\bar a, \bar b)$.
\item
 $T\vdash \forall \bar x\,\bar y\, 
 ((\psi(\bar x)\wedge\theta_1(\bar x, \bar y))\, 
 \rightarrow
  \theta_2(\bar x, \bar y))$.
 \end{enumerate}
\end{thm}

\preuve
Let $A$ be a pac model of\  $T$, $\bar a\in A$, and $\psi$  a positive
formula such that $A\nvDash\psi(\bar a)$.  
 Since every embedding from $A$ into 
a model of\  $T$ is an immersion, then the  set  
  of h-inductive sentences $\{T, Diag(A), \psi(\bar a)\}$
  is inconsistent. 
 By compactness,  there exist 
$\theta_1(\bar a, \bar b), \neg\theta_2(\bar a, \bar b)\in Diag(A)$
such that 
$\{T,\theta_1(\bar a, \bar b),
 \neg\theta_2(\bar a, \bar b), \psi(\bar a)\}$
is inconsistent, thereby we obtain  
$$T\vdash \forall \bar x\bar y\,
 ((\psi(\bar x)\wedge\theta_1(\bar x, \bar y))\, 
 \rightarrow
  \theta_2(\bar x, \bar y)).$$
  
For the other direction, let $A$ be a model of\  $T$ that satisfies 
the properties $(1)$ and $(2)$ of the theorem.
 Let 
 $p$ be  an embedding from $A$ into a model $B$ of\  $T$.  
 Let $\psi$ be a 
positive formula and $\bar a\in A$ such that $A\nvDash \psi(\bar a)$.
By hypothesis, there exist $\theta_1, \theta_2$  a pair of  positive quantifier-free formulas and  $\bar b\in A$
that satisfy the following properties: 
$$
 \left\{
    \begin{array}{ll}
        A\models \theta_1(\bar a, \bar b) \\
       A \nvDash\theta_2(\bar a, \bar b) \\
        T\vdash \forall \bar x\bar y\, 
 ((\psi(\bar x)\wedge\theta_1(\bar x, \bar y))\, 
 \rightarrow\theta_2(\bar x, \bar y))
    \end{array}
\right.
$$
Now since  $p$ is an embedding  then
 $B\models \theta_1(\bar a, \bar b)$
and $B\nvDash\theta_2(\bar a, \bar b)$. Consequently  
 $B\nvDash \psi(\bar a)$,  so $p$ is an immersion.\qed
  
\begin{df}\label{weakresultant}
Let $T$ be an h-inductive theory and $\psi$  a positive 
formula. We denote by $Alc_T(\psi)$  the set 
of  pairs of positive quantifier-free formulas
$(\theta_1, \theta_2)$  that satisfy the following property:\\
 there 
exist  a pac model $A$ of\  $T$,  $\bar a$ and $\bar b$ 
tuples  from $A$ such that:
$$
(\star) \left\{
    \begin{array}{ll}
    A\nvDash\psi(\bar a)\\
        A\models \theta_1(\bar a, \bar b) \\
       A \nvDash\theta_2(\bar a, \bar b) \\
        T\vdash \forall \bar x\bar y\, 
 ((\psi(\bar x)\wedge\theta_1(\bar x, \bar y))\, 
 \rightarrow\theta_2(\bar x, \bar y))
    \end{array}
\right.
$$

\begin{itemize}
\item Let $\Delta$ be a subset of\  $Alc_T(\psi)$. 
We say that 
  $Alc_T(\psi)$ 
  is equivalent to $\Delta$ and we write 
  $Alc_T(\psi)\equiv \Delta$  if for every 
  pac model $A$ of\  $T$  and $\bar a\in A$; 
  if $A\nvDash\psi(\bar a)$
  then there are   $(\theta_1, \theta_2)\in\Delta$
   and $\bar b\in A$,
   that satisfy the property $(\star)$ above.
   \item Let $\Delta$ be a subset of\  
  $Ctr_T(\psi)$. We say that 
   $Ctr_T(\psi)$
  is equivalent to $\Delta$ and we write 
  $Ctr_T(\psi)\equiv \Delta$  if  for every 
   pc model $A$ of\  $T$ and $\bar a\in A$; 
  if $A\nvDash\psi(\bar a)$
  then there are
   $\phi\in\Delta$ and $\bar b\in A$ such that 
    $A\models\phi(\bar a)$.
\end{itemize}
  
\end{df}
\begin{rem}\label{disjCag}
\begin{itemize}
\item 
  If $ (\alpha_1, \alpha_2), (\beta_1, \beta_2)\in 
 Alc_T(\psi)$, we have 
 $(\alpha_1\vee\beta_1, \alpha_2\vee\beta_2)\in Alc_T(\psi)$. 
 \item An h-inductive theory $T$ has a model-companion if and 
 only if for every positive formula $\psi$ there exists 
 $\psi^c\in Ctr_T(\psi)$ such that $Ctr_T(\psi)\equiv\{\varphi\}$.\\
 In this case the class of pc  models is axiomatized by $T_k(T)$
 the following h-inductive theory:
 $$T\cup\{\forall\bar x\, (\psi(\bar x)\vee \psi^c(\bar x))
 |\ \psi\ \text{ranges over the set of positive formulas}\}.$$
\end{itemize}
\end{rem}
\begin{df}
Let $T$ be an h-inductive theory. We say that $T$ has an algebraic 
model-companion 
if the class of pac models of\  $T$
 is axiomatized by a first order  
theory.
\end{df}
\begin{thm}\label{axiomweakpec}
An h-inductive theory $T$  has an algebraic 
model-companion 
 if and only if for every positive 
formula $\psi$,  $Alc_T(\psi)$ is equivalent to a pair
$(\phi_1, \phi_2)$ of\  $Alc_T(\psi)$. In this case  the class 
of pac  models of\  $T$ is axiomatized by $T_h$ the following 
inductive theory: 
$$T\cup\{\forall \bar x\,\bar y\  
 ((\psi(\bar x)\wedge\phi_1(\bar x, \bar y))\, 
 \rightarrow
  \phi_2(\bar x, \bar y)),\,\forall\bar x\, (\psi(\bar x)
  \vee
  (\exists\bar y\, (\phi_1(\bar x, \bar y)\wedge\neg\phi_2(\bar x,\bar y))))
  \}$$ where $\psi$ ranges over the set of positive formulas. 
\end{thm}
\preuve 
Suppose that  the class of pac models of\  $T$
is axiomatized by a first order 
theory $T^\star$. Let 
$\psi$ be a positive formula such that  $Alc_T(\psi)$ is not
 equivalent to any of its 
finite subset. Then, for  every finite subset of  $Alc_T(\psi)$ 
(which by the remark \ref{disjCag} can be represented by a pair 
$(\theta_1, \theta_2)\in Alc_T(\psi)$) there are  a model
 $A$ of\  $T^\star$ and $\bar a\in A$, such that
 $$
\left\{
    \begin{array}{ll}
        A\nvDash\psi(\bar a) \\
        A\models\exists\bar x(\theta_1(\bar a, \bar x)\wedge\theta_2(\bar a, \bar x)).
\end{array}
\right.
$$ 
By compactness, there exist $B$ a model of\  $T^\star$ and
 $\bar b\in B$
such that for every $(\theta_1, \theta_2)\in Alc_T(\psi)$ we have:
$$
\left\{
    \begin{array}{ll}
        B\nvDash\psi(\bar a) \\
      B\models\exists\bar x(\theta_1(\bar a, \bar x)
      \wedge\theta_2(\bar a, \bar x)).
\end{array}
\right.
$$  
Contradiction from the syntactic characterization of 
pac models. Thereby,  
there is $(\theta_1, \theta_2)\in Alc_T(\psi)$ such that 
$Alc_T(\psi)\equiv\{(\theta_1, \theta_2)\}$. It is clear that 
$T^\star$ is logically equivalent to the theory $T_h$ given in the 
theorem.

For the other direction,  
suppose that for each positive formula $\psi$ there 
is $(\phi_1, \phi_2)\in Alc_T(\psi)$ such that 
$Alc_T(\psi)\equiv\{(\phi_1, \phi_2)\}$. Let $T_h$ be the 
theory $$T\cup\{\forall \bar x\,\bar y 
 ((\psi(\bar x)\wedge\phi_1(\bar x, \bar y))\, 
 \rightarrow
  \phi_2(\bar x, \bar y)),\,\forall\bar x\, (\psi(\bar x)
  \vee
  (\exists\bar y\, \phi_1(\bar x, \bar y)\wedge\neg\phi_2(\bar x,\bar y)))
  \}$$ where $\psi$ runs through the set of positive formulas.\\
  It is clear that every pac model of\  $T$ is a model of\  $T_h$.
  Conversely,  
  let $A$ be a model of\  $T_h$ and $e$ an embedding from $A$ into 
  $B$ a model of\  $T$. Let $\bar a\in B$ and $\psi $ be a positive 
  positive formula such that  $A\nvDash \psi(\bar a)$. By
  hypothesis, there exist 
  $(\phi_1, \phi_2)$ a pair of positive quantifier-free formulas  and 
  $\bar b\in A$, such that $A\models \phi_1(\bar a,\bar b)$ and 
  $A\nvDash \phi_2(\bar a,\bar b)$. Given that $e$ is an embedding 
  it follows that  $B\models \phi_1(\bar a,\bar b)$, and 
  $B\nvDash \phi_2(\bar a,\bar b)$. Since   
  $T\vdash\forall \bar x\,\bar y 
 ((\psi(\bar x)\wedge\phi_1(\bar x, \bar y))\, 
 \rightarrow
  \phi_2(\bar x, \bar y))$ we obtain 
  $B\nvDash\psi(\bar a)$, thus $e$ is an immersion and 
   $A$ is a pac model of\  $T$. Thereby the class of pac models of\  $T$
   is axiomatized by the theory $T_h$.\qed
   \begin{rem}
  If the class of pac models of an h-inductive theory is axiomatized 
  by an h-inductive theory, then every model of $T_k(T)$ is 
  a pac model of $T$. 
   \end{rem}

 \begin{lem}\label{lemma 9}
Let $A$ be a model of\  $T$ which is immersed in a 
pac model $C$  of\  $T$. If $A$  has
the amalgamation property for embeddings 
then $A$ is a pac model of\  $T$.
\end{lem}
\preuve Let $e$ be an embedding from $A$ into $B$ a model of\  $T$.
since $A$ has the amalgamation of embeddings, we obtain the following
commutative diagram:
\[
\xymatrix{
    A \ar[r]^{e} \ar[d]_{i_m} & {B} \ar[d]^{e_2} \\
    C \ar[r]_{e_3} & {D}
  }
\]  
where $D$ is a model of\  $T$,  $e_2$ and  $e_3$  embeddings.
 Given that  
$C$ is a pac model  of\  $T$ and $D$  a model of\  $T$, then $e_3$
is a immersion. Since  $e_3\circ i_m= e_2\circ e$ and $e_3\circ i_m$
is an immersion, then $ e_2\circ e$ is an immersion, so 
$e$ is an immersion.   Thus $A$ is a pac model of\  $T$.\qed
\begin{rem}
Note that if $T$ is an inductive theory  in Robinson's sense, and 
 a structure $A$ immersed in some model of $T$, then 
$A$ is not necessary a model of $T$. Given that the 
the characterization of pac models of an h-inductive
 theory is made by inductive sentences.
 So it is natural  to ask the question whether  
 the condition of  amalgamation 
 of embeddings in the lemma \ref{lemma 9} is a 
 necessary and sufficient condition. 

\end{rem}
\subsection{E-elementary  extension}
In this section, we introduce and make preliminary study of the notion of e-elementary  extension  inspired by the the notion of elementary 
extension given by Poizat in \cite{poizat}.\\ 
Recall that $B$ is an elementary extension of $B$ in the 
the sense of Poizat if and only if $B$ is a pc model 
of $T_i(A)$.

Let $A$ be a $L$-structure. We use the notation $T_A$ for 
the h-inductive $L(A)$-theory $\{T_u(A), Diag^+(A)\}$.\\ 

 \begin{lem}\label{intralgebextens}
 Let $A$ be a $L$-structure. The theories  
$T_A $  and $T_i(A)$ are e-companions
$L(A)$-theories.
\end{lem}
\preuve
 Let $B_1$ be a pac model of\  $T_A$ and 
$B_2$  a pac model of\  $T_i(A)$. 
 Given that $B_2$ is a model of\  $T_i(A)$ and 
$A$ is immersed in both $B_1$ and $B_2$, then by the lemma
\ref{pamalgam}, we obtain the following diagram
\[
\xymatrix{
    A \ar[r]^{i_m} \ar[d]_{i_m} & {B_2} \ar[d]^{i} \\
    {B_1} \ar[r]_{e} & {D}
  }
\]  
where $D$ is a model of\  $T_i(B_2)$. Therefore 
$B_1$ is embedded in a model of\  $T_i(A)$, and $B_2$ is embedded in a 
model of\  $T_A$. \qed
\begin{df}\label{e-ex}
Let $A$ and $B$ be two $L$-structures.  $B$ is said to be  an 
e-elementary
extension of\  $A$, in symbols $A\prec_eB$, if  $B$ is a pac model of\  $T_A$.
\end{df}
\begin{rem}
From the lemma \ref{intralgebextens} and the definition \ref{e-ex} 
it follows that 
$A\prec_e B$ if and only if $B$ is a pac model of $T_i(A)$. Since 
every pc model of  an h-inductive theory is a  pac model, it follows that the notion  
of e-elementary extension is weaker than  the notion of elementary extension
given by Poizat.  
\end{rem}
The following lemma gives some  properties of the  notion
of e-elementary extension. 
\begin{lem}\label{propeleext}
\begin{enumerate}
\item For every $L$-structure $A$ we have $A\prec_eA$.
\item Let $T$ be an h-inductive theory, $A$ a pac model of\  $T$, and 
$B$ an elementary extension of\  $A$. 
Then $B$ is a pac model of\  $T$.
\item Let $A, B, C$ three L-structures such that, 
$A\prec_eB\prec_eC$. Then $A\prec_eC$.
\item If $A\prec_eB$ then $T_i(A)= T_i(B)$ in the language $L(A)$.
\item If $A\prec_eC$, $B\prec_eC$ and $A$ is embedded in $B$. Then 
$A\prec_eB$.
\end{enumerate}
 \end{lem}
\preuve 
\begin{enumerate}
\item Clear.
\item Consider the following diagram 
$$A\xrightarrow{i_m}B\xrightarrow{p}C,$$ 
where  $C$ is a model of\  $T$ and $p$ an embedding. In order to prove 
that $B$ is a pac model of\  $T$ it suffices to show that 
$p$ is an immersion. \\
Given that $A$ is a pac model of\  $T$ then $p\circ i_m$ is an immersion, so $C$ is 
a model of\  $T_A$, which implies that $p$ is an immersion.
\item Suppose that $C$ is embedded in a model $D$ of\  $T_A$. We 
have the following diagram:
$$A\xrightarrow{i_1}B\xrightarrow{i_2}C\xrightarrow{p}D$$
where $p$ is an embedding, $i_1$ and  $i_2$ are immersions.
Since $A\prec_eB$, $p\circ i_2$ is am immersion, so $D$ is 
a model of\  $T_B$. As $B\prec_eC$ then $p$ is an immersion.
\item Since $A$ is immersed in $B$ then  $T_i(B)\subseteq T_i(A)$
over the language $L(A)$. On the other hand, as 
$T_A$ and $T_i(A)$ are  e-companions theories and 
$A\prec_eB$ then $B\models T_i(A)$, which implies 
$T_i(A)\subseteq T_i(B)$.
\item According to the hypothesis of the proposition ($5$) we 
have the following diagram:
$$A\xrightarrow{p}B\xrightarrow{i}C$$
where $p$ is an embedding and $i$ an immersion.\\
Given that $A\prec_e C$, then $i\circ p$ is an immersion, so 
$p$ is an immersion. By the proposition ($4$) we have 
$T_i(A)=T_i(C)=T_i(B)$ over the language $L(A)$. Since 
$C$ is a pac model of\  $T_A$, $B$ is immersed in $C$ and $C$ is a 
model of\  $T_i(B)$ in the language $L(A)$, if follows from the lemma 
\ref{substrextension} that $B$ is a  pac model of\  $T_A$.\qed
\end{enumerate}

\section*{References}
\bibliographystyle{plain}
\bibliography{reference}

\begin{thebibliography}{1}

\bibitem{ana}
Mohammed Belkasmi.
\newblock Positive model theory and amalgamations.
\newblock {\em Notre Dame Journal of Formal Logic}, 55(2), 2014.

\bibitem{Almaz}
Almaz Kungozhin.
\newblock Existentially closed and maximal models in positive logic.
\newblock {\em Algebra and Logic}, 51(6), 2013.

\bibitem{poizat}
Bruno Poizat.
\newblock Univers positifs.
\newblock {\em Journal of Symbolic Logic}, 71(3), 2006.

\bibitem{rob1}
Abraham Robinson.
\newblock On the notion of algebraic closedness for noncommutative groups and
  fields.
\newblock {\em Journal of Symbolic Logic}, 36:441--444, 1971.

\bibitem{scott}
W.R. Scott.
\newblock Algebraically closed groups.
\newblock {\em Proceedings AMS}, 2(1), 1951.

\bibitem{begnacpoizat}
Itai~Ben Yaacov and Bruno Poizat.
\newblock Fondements de la logique positive.
\newblock {\em Symbolic Logic}, 72(4):1141--1162, 2007.

\end{thebibliography}

\end{document}